\documentclass{amsart} 
\usepackage{amssymb,latexsym}

\numberwithin{equation}{section}

\newtheorem{theorem}{Theorem}[section]
\newtheorem{proposition}[theorem]{Proposition}
\newtheorem{corollary}[theorem]{Corollary}
\newtheorem{remark}[theorem]{Remark}
\newtheorem{lemma}[theorem]{Lemma}
\newtheorem{example}[theorem]{Example}
\newtheorem{definition}[theorem]{Definition}

\def\proof{\smallskip\noindent {\bf Proof. }}
\def\endproof{\hfill$\square$\medskip}

\newcommand{\journal}[1]{{\sl #1}}
\newcommand{\arttitle}[1]{{\rm #1}}

\def\T{\mathcal{T}}
\def\k{{k}}
\def\l{\ell}
\def\wnot{w_\mathrm{o}}

\begin{document}



\title{
Mixed Bruhat operators
and\\Yang-Baxter equations 
for Weyl groups
}

\author{Francesco Brenti}
\address{\noindent Dipartimento di Matematica, 
Universit{\'a} di Roma ``Tor Vergata'', 
Via della Ri\-cer\-ca Scientifica, 00133 Roma, Italy} 
\email{brenti@mat.uniroma2.it} 

\author{Sergey Fomin}
\address{Department of Mathematics, Massachusetts Institute of
  Technology, Cambridge, Massachusetts 02139}
\email{fomin@math.mit.edu}

\author{Alexander Postnikov}
\address{Department of Mathematics, Massachusetts Institute of
  Technology, Cambridge, Massachusetts 02139}
\email{apost@math.mit.edu}

\thanks{
The authors were supported in part by MSRI (NSF grant \#DMS-9022140).
The second author was also supported by NSF grant \#DMS-9700927.}

\date{April 8, 1998
}



\maketitle

\section{Introduction}
\label{sec:intro}

In this paper, we introduce and study a family of 
operators which act in the span of a Weyl group~$W$ 
and provide a multi-parameter solution
to the quantum Yang-Baxter equations of the corresponding type. 
These operators are then used to obtain new results (as well as new
proofs of the known facts) concerning the Bruhat order of~$W$. 

Let us briefly review the contents of the paper.
Section~\ref{sec:ybe} is devoted to preliminaries related to
Coxeter groups and associated Yang-Baxter equations. 
The \emph{mixed Bruhat operators} $M_\tau$ 
are defined in Section~\ref{sec:mbo} by the formula 
\[
M_\tau(w) = \left\{ 
\begin{array}{ll}
p_\tau \tau w  & \mbox{if\ $\l(\tau w)>\l(w)$;} \\[.07in]
q_\tau \tau w   & \mbox{if\ $\l(\tau w)<\l(w)$,} \\
\end{array} \right.
\]
where $p_\tau$ and $q_\tau$
are scalar parameters that depend on a reflection~$\tau$.
In Section~\ref{sec:mbo}, we also state our main result, which
describes a family of choices for the parameters $p_\tau$
and $q_\tau$ 
such that the associated operators 
$R_\tau=1+\varepsilon M_\tau$ satisfy the Yang-Baxter equations. 
The well known $W$-analogue of the
classical Yang's solution 
can be obtained from our
solution as a particular degeneration. 

In Section~\ref{sec:Q}, we define the \emph{quantum Bruhat operators}
as certain limiting case of 
the~$M_\tau\,$.  
These operators 
play an important role (which we briefly 
explain) in the explicit description of the multiplicative structure
of the (small) quantum cohomology ring of the flag manifold of the
corresponding type. 

Section~\ref{sec:proof} contains the proof of the main result.
To make the presentation more accessible, we first treat the case of
the symmetric group, and then provide the proof in complete
generality. 
We also show how our family of solutions arises naturally in an
attempt to solve the Yang-Baxter equations within this class of
operators. 

Section~\ref{sec:combin} is devoted to combinatorial applications of
our operators. 
For an arbitrary element $u\in W$, we define 
a graded poset with unique minimal element~$u$ that we call a
\emph{tilted Bruhat order}. 
The usual Bruhat order corresponds to the special case~$u=e$. 
We then prove that tilted Bruhat orders are lexicographically
shellable graded posets whose 
every interval is Eulerian. 
This generalizes well known results of Verma, Bj\"orner, Wachs, and
Dyer.

\section{Coxeter groups and Yang-Baxter equations}
\label{sec:ybe}

We first introduce some standard terminology and notation related to
Coxeter groups and root systems. 
In what follows, $W$ is a \emph{Weyl} group, and $S$ is the set of its
\emph{simple reflections}. 
In particular, $(W,S)$ is a \emph{Coxeter system}, i.e, the 
finite (Coxeter) group $W$ is given by the presentation 
\[
(st)^{m(s,t)}=1\,, \ \ s,t\in S\,, 
\]
where the nonnegative integers $m(s,t)$ satisfy $m(s,s)=1$
and $m(s,t)=m(t,s)>1$ for $s\neq t$.
The group $W$ is finite; let $\wnot$ 
denote its longest element.
Most constructions in this section can be extended naturally 
to arbitrary Coxeter groups.

For an element $w\in W$, an expansion 
$w=s_1\cdots s_l $
of minimal possible length $l$ is called a \emph{reduced
decomposition}. The number $l=\l(w)$ is the \emph{length} of~$w$. 
The elements of the set 
$T =\{ w s w^{-1} : w \in W, \; s \in S \}$ are the \emph{reflections}
of~$W$.

The \emph{Bruhat order} on $W$ is defined as follows: 
$u \leq v$ if and only if there exist 
$t_{1},\ldots , t_{r} \in T$ such that 
$t_{r} \ldots  t_{1}\,u=v$ and 
$\l(t_{i} \ldots t_{1}\,u ) > \l(t_{i-1} \ldots t_{1}\,u)$ for $i=1,
\ldots ,r$. 

Geometrically, the group $W$ can be represented in terms of a 
\emph{root system}~$\Phi$.
A subgroup $W'$ of $W$ generated by a subset $A \subseteq T$
is called a \emph{reflection subgroup};
it corresponds to a root subsystem of~$\Phi$. 
Hence  $W'$ is again a Coxeter group, with the set of \emph{canonical}
(Coxeter) \emph{generators} $S'$ corresponding to the simple roots of this
subsystem. (This can be extended to any Coxeter group: 
see~\cite{Deo2,Dy2} or \cite[8.2]{Hum}.)  
We will only be interested in the case where 
$W'$ is a \emph{dihedral reflection subgroup},
i.e., $S'$ has two elements. 
A dihedral reflection subgroup is \emph{maximal} if it is not
contained in another such subgroup.
Maximal dihedral subgroups correspond
to two-dimensional root subsystems obtained by intersecting $\Phi$
with a plane spanned by a pair of positive roots.

Let $N=\l(\wnot)$.
Following Dyer~\cite{Dy}, we say that a bijection 
$\varphi : T \rightarrow \{1,\dots,N\}$ is a (total) 
\emph{reflection ordering} if, for any dihedral reflection
subgroup $W'$ with canonical generators $a$ and~$b$, 
the sequence 
\[
\varphi (a) , \varphi(aba) , \varphi (ababa) , \ldots 
,\varphi (babab), \varphi (bab),\varphi (b)
\]
is either increasing or decreasing. 
(It is enough to require this for every maximal dihedral subgroup.) 
Reflection orderings correspond (bijectively) to reduced
decompositions of $\wnot$ in the following standard way: 
$\varphi $ is a reflection  ordering if and only if there exists  a 
reduced decomposition 
$\wnot=s_1 \ldots s_N$ such that
\begin{equation} 
\varphi  ^{-1} (j)=s_N s_{N-1} \ldots 
s_{j+1} s_j s_{j+1} \ldots s_{N-1} s_{N} 
\label{eq:ref-ord}
\end{equation}
for $j=1,\ldots N$.

\begin{definition}
\label{def:ybe}{\rm
A family $\{R_\tau\}_{\tau\in T}$ of elements of a monoid 
is called an (extensible) \emph{solution to the Yang-Baxter equations}
for~$W$ if for any dihedral reflection subgroup~$W'$ of~$W$ 
with canonical generators $a$ and $b$, we have
\begin{equation}
\label{eq:ybe}
R_a  R_{aba}  R_{ababa} \cdots R_{bab} 
R_b = R_b  R_{bab} \cdots R_{ababa} 
R_{aba}  R _a \,.
\end{equation}
In particular, 
if $a,b \in T$ and $ab=ba$, then $R_a R_b =R_b R_a\,$. 
The collection $\{ R_t \}_{t \in T}$ satisfying
the Yang-Baxter equations~(\ref{eq:ybe}) is frequently called an
(extensible) $R$-\emph{matrix} (of the corresponding type); 
we will not use this terminology here.
}\end{definition}

The definition above makes sense for any finite Coxeter group.   
In the case of a Weyl group, equations~(\ref{eq:ybe}), 
stated case by case in terms of 
the root system for~$W$, were given by Cherednik
(implicit in~\cite{cherednik1} and explicit
in~\cite[Definition~2.1a]{cherednik2}),
along with a number of solutions. 

\begin{remark}{\rm
The word ``extensible'' (which we will later omit;
cf.\ \cite[Definition~2.2]{cherednik2}) 
indicates that we ask for~(\ref{eq:ybe}) to be satisfied for
\emph{all} dihedral subgroups, not just for the \emph{maximal} ones. 
(The distinction is only relevant in non-simply-laced cases.) 
This stronger condition, however not needed for the general ``Yang-Baxter
machinery'' to work, will actually be satisfied by all
solutions constructed in this paper, which explains our choice of 
definition.
}\end{remark}

For the type $A_{n-1}\,$, the Weyl group is the symmetric
group~$S_n\,$, the set $T$ consists of all transpositions~$(ij)\in
S_n\,$, and the equations (\ref{eq:ybe})
are the celebrated (quantum) 
Yang-Baxter equations (see, e.g.,~\cite{ybe}). Let us explain. 
Let $R_{ij}$ be a shorthand for~$R_{(ij)}\,$.
Then (\ref{eq:ybe}) becomes  
\begin{equation}
\qquad R_{ij}R_{kl} = R_{kl}R_{ij} 
\qquad \textrm{ if $i,j,k,l$ are distinct;}
\label{eq:ybe-a1}
\end{equation}
\begin{equation}
R_{ij}R_{ik}R_{jk} = R_{jk}R_{ik}R_{ij}\qquad 
\textrm{ if $i<j<k$.}\qquad  
\label{eq:ybe-a2}
\end{equation}

\begin{example}{\rm
The first solution to the Yang-Baxter equations
was given by Yang 
in his pioneering paper~\cite{yang}, 
where he observed that the elements
\begin{equation}
\label{eq:yang}
R_{ij}=1+\frac{(ij)}{x_j-x_i}
\end{equation}
of the group algebra of the symmetric group~$S_n$ 
satisfy~(\ref{eq:ybe-a1})--(\ref{eq:ybe-a2}),
for any choice of distinct parameters~$x_1,\dots,x_n\,$.
This generalizes to an arbitrary Weyl group as
follows~\cite{cherednik1,cherednik2}: 
\begin{equation}
\label{eq:yang-general}
R_\tau=1+\frac{\varkappa_\tau \tau}{\langle x,\alpha\rangle }\, ,
\end{equation}
where $\alpha$ is the positive root corresponding to~$\tau$,
and $\varkappa_\tau$ is a scalar whose value only depends on whether the
root~$\alpha$ is short or long. 
(In other words, $\varkappa_\tau=\varkappa_\sigma$ if reflections
$\tau$ and $\sigma$ are conjugate to each other.)
}\end{example}

The fact that every two reduced decompositions 
of the element $\wnot\in W$ are related 
by a sequence of elementary Coxeter transformations
(see, e.g., ~\cite[Section~8.1]{Hum})
translates (using~(\ref{eq:ref-ord})) into every two reflection
orderings being related by a sequence of Yang-Baxter-type moves
of the form
\[
\dots,a,aba,ababa,\dots,bab,b,\dots\ \leadsto\ 
\dots,b,bab,\dots,ababa,aba,a,\dots\,,
\]
where $a$ and $b$ are the canonical generators of the (maximal)
dihedral subgroup they generate. 
This implies the following statement. 

\begin{proposition}
\label{pro:invariant-product}
Let $\{R_\tau\}$ be a solution of the Yang-Baxter equations 
for a finite Coxeter group~$W$,
and let $\varphi:T\to \{1,\dots,N\}$ be a reflection ordering on~$T$. 
Then the product
\begin{equation*}
\label{eq:R-product} 
\prod _{i=1}^{N} R_{\varphi ^{-1} (i)}
= R_{\varphi ^{-1} (1)} \cdots  R_{\varphi ^{-1} (\l(\wnot)}  
\end{equation*}
does not depend on the choice of a reflection ordering~$\varphi$.
\end{proposition}

\section{Mixed Bruhat operators}
\label{sec:mbo}

We will work over a ground field~$\k$ of characteristic~0. 

\begin{definition}
\label{def:mbo}{\rm
Let $\{p_\tau\}$ and $\{q_\tau\}$
be two families of scalar parameters indexed by reflections $\tau\in T$.
The \emph{mixed Bruhat operators} $M_\tau$ are 
 linear operators acting in the $\k$-span  $\k[W]$ of the group~$W$ by
\begin{equation}
\label{eq:M}
M_\tau(w) = \left\{ 
\begin{array}{ll}
p_\tau \tau w  & \mbox{if\ $\l(\tau w)>\l(w)$;} \\[.1in]
q_\tau \tau w   & \mbox{if\ $\l(\tau w)<\l(w)$.} \\
\end{array} \right.
\end{equation}
Let $\varepsilon$ be a formal variable with values in~$\k$,
and define the operators
\begin{equation}
\label{eq:R}
R_\tau = 1 + \varepsilon M_\tau\ .
\end{equation} 
}\end{definition}

We will now describe a particular multi-parametric construction that
allows to choose the $p_\tau$ and  $q_\tau$ 
so that the operators $R_\tau$ 
satisfy the Yang-Baxter equations~(\ref{eq:ybe}). 

\begin{definition}{\rm
A function $\alpha\mapsto E(\alpha)$ defined on the set of positive
roots is called \emph{multiplicative} if, whenever $\alpha$, $\beta$,
and $\alpha+\beta$ are positive roots, we have
\begin{equation}
\label{eq:multiplicativity}
E(\alpha+\beta)=E(\alpha)\,E(\beta)\,.
\end{equation}
To construct such a function, simply assign arbitrary values to the
simple roots, and then extend by multiplicativity. 
A typical example of a multiplicative function is given by 
\begin{equation}
\label{eq:exponential} 
E(\alpha)=e^{\langle\alpha,x\rangle}\,,
\end{equation}
where $x$ is an arbitrary vector. 
Notice, however, that (\ref{eq:exponential}) does not allow for
$E(\alpha)=0$, a possibility that we do not want to exclude. 
}
\end{definition}

\begin{theorem}
\label{th:main-W}
Let $\alpha\mapsto E_1(\alpha)$ and $\alpha\mapsto E_2(\alpha)$ be 
multiplicative functions on the set
of positive roots such that $E_1(\alpha)\neq E_2(\alpha)$ for
every~$\alpha$. 
Let $\kappa_\tau$ be a scalar whose value only depends on the length
of the positive root~$\alpha$ corresponding to~$\tau$. 
Define parameters 
$p_\tau\,,q_\tau\,$, for $\tau\in T$, by
\begin{equation}
\label{eq:global-p}
p_\tau = \frac{\kappa_\tau E_1(\alpha)}{E_1(\alpha)-E_2(\alpha)}  
\end{equation}
and 
\begin{equation}
\label{eq:global-q}
q_\tau = \frac{\kappa_\tau E_2(\alpha)}{E_1(\alpha)-E_2(\alpha)} \,.
\end{equation}
Then the operators $R_\tau$ given by {\rm (\ref{eq:M})--(\ref{eq:R})}
satisfy the quantum Yang-Baxter equations~{\rm (\ref{eq:ybe})}. 
\end{theorem}

\begin{example}{\rm 
Consider the type $A_{n-1}$ case where $W$ is the symmetric
group~$S_n\,$. 
For a reflection $\tau=(ij)$, 
we will use the notation $p_{ij}\,$, $q_{ij}\,$, 
$M_{ij}\,$, and $R_{ij}$ 
instead of $p_\tau\,$, $q_\tau\,$, $M_\tau\,$, and~$R_\tau\,$. 
Hence 
\begin{equation}
\label{eq:RA}
R_{ij}(w) = \left\{ 
\begin{array}{ll}
w +\varepsilon p_{ij} \tau w  & \mbox{if\ $\l(\tau w)>\l(w)$;} \\[.1in]
w +\varepsilon q_{ij} \tau w   & \mbox{if\ $\l(\tau w)<\l(w)$.} \\
\end{array} \right.
\end{equation}

The positive root corresponding to $\tau=(ij)$ is 
$\alpha=\alpha_i+\cdots+\alpha_{j-1}\,$,
where $\alpha_1,\dots,\alpha_{n-1}$ are the simple roots, ordered in a
standard way. 
Thus the multiplicative functions $E_1$
and $E_2$ are determined by the values
$p_i=E_1(\alpha_i)$ and $q_i=E_2(\alpha_i)$, as follows: 
$E_1(\alpha)=p_i\cdots p_{j-1}\,$, 
$E_2(\alpha)=q_i\cdots q_{j-1}\,$.
This leads to 
\begin{equation}
\label{eq:global-pA}
p_{ij} = \frac{\kappa\, p_i\cdots p_{j-1}}{p_i\cdots p_{j-1}-q_i\cdots
  q_{j-1}}  
\end{equation}
and 
\begin{equation}
\label{eq:global-qA}
q_{ij} = \frac{\kappa\, q_i\cdots q_{j-1}}
{p_i\cdots p_{j-1}-q_i\cdots q_{j-1}} \,.
\end{equation}
(Since all roots have the same length, we drop the subscript~$\tau$
in~$\kappa_\tau\,$.) 
Substituting this into (\ref{eq:RA}), we obtain a family of solutions
of the Yang-Baxter equations (of type~$A$). 
In formulas  (\ref{eq:global-pA})--(\ref{eq:global-qA}), 
$\kappa$ is an arbitrary scalar,
while the parameters $p_i$ and $q_i$ should be chosen so that none of the
denominators vanish.
Notice that we do not use a single set of parameters $t_i=q_i/p_i$ in
order to, first, keep the symmetry between the $p_i$ and the $q_i$ 
and, second, allow for the possibility of $p_i=0$. 
}\end{example}

\begin{remark}{\rm
The analogue (\ref{eq:yang-general})
of Yang's solution of the Yang-Baxter equation
can be obtained from the solution given in
Theorem~\ref{th:main-W} as a particular limiting case. 
Let $\kappa_\tau=\delta\,\varkappa_\tau\,$, where $\delta$ is a scalar. 
Fix a vector~$x$,
and set $E_1(\alpha)=e^{\delta\langle\alpha,x\rangle}$ 
and $E_2(\alpha)=1$. 
Making these substitutions into
(\ref{eq:global-p})--(\ref{eq:global-q}) and taking the limit as
$\delta\to 0$, we obtain
$
p_\tau = q_\tau = 
\displaystyle\frac{\varkappa_\tau}{\langle\alpha,x\rangle} \,, 
$
which means that the operators $\lim_{\delta\to 0} M_\tau$ act by left
multiplication by 
$\displaystyle\frac{\varkappa_\tau \tau}{\langle\alpha,x\rangle}$,
as desired.  
}\end{remark}

\section{Rescaling. Quantum Bruhat operators}
\label{sec:Q}

Rescaling is a very simple yet sometimes helpful way of producing 
new solutions  to the Yang-Baxter equations from existing ones. 
In this section, we show how rescaling of the mixed Bruhat operators
leads in the limiting case to the construction of 
``quantum Bruhat operators'' for an arbitrary Weyl group~$W$. 
These operators, introduced in~\cite{FK} for type~$A$, 
appear in the analogue of Monk's formula for the (small) quantum
cohomology ring of the flag manifold (see below). 
In this paper, we are mainly concerned with their combinatorial
applications. 

Suppose that $\{M_\tau\}_{\tau\in T}$ is a family of mixed Bruhat
operators such that the corresponding  operators 
$R_\tau=1+\varepsilon M_\tau$ satisfy 
the Yang-Baxter equations~(\ref{eq:ybe}). 
Let $\{\gamma_w\,:\,w\in W\}$ be a collection of nonzero scalars. 
Then the \emph{rescaled} operators $\widetilde M_\tau$ defined by
\begin{eqnarray}
\label{eq:Mtilde}
\begin{array}{l}
\widetilde M_\tau(w)
=\displaystyle\frac{\gamma_{\tau w}}{\gamma_w} \, M_\tau(w)
\end{array}
\end{eqnarray}
are also a solution to (\ref{eq:ybe}).
This follows from the fact that 
$M_\tau(w)$ is always a scalar multiple of~$\tau w$,
and therefore 
$\widetilde M_\tau=\Gamma M_\tau \Gamma^{-1}$,
where $\Gamma(v)=\gamma_v v$ for $v\in W$. 

Let ${\rm ht}(\alpha)$ denote the \emph{height} of a positive
root~$\alpha$, i.e., the sum of the coefficients in the expansion
of~$\alpha$ in the basis of simple roots. 
Then for any scalar~$h$ and any multiplicative function $\alpha\mapsto
E(\alpha)$, the function $\alpha\mapsto h^{{\rm ht}(\alpha)}E(\alpha)$
  is also multiplicative.

Let $\delta\neq 0$ be a scalar parameter (eventually, we will take
$\delta\to 0$),
and let $\alpha\mapsto E(\alpha)$ be a multiplicative function. 
Let the parameters $p_\tau$ and $q_\tau$ of the mixed Bruhat operators
$M_\tau$ be given by (\ref{eq:global-p})--(\ref{eq:global-q}) with 
\begin{eqnarray}
\label{eq:specialize}
\begin{array}{rcll}
\kappa_\tau&=&\delta^{-1}\ ,\\[.1in]
E_1(\alpha) &=& 1 \ ,\\[.1in]
E_2(\alpha) &=& \delta^{2{\rm ht}(\alpha)}E(\alpha)\ .
\end{array}
\end{eqnarray}
Using notation $F\approx G$ for 
$\displaystyle\lim_{\delta\to 0}{F/G}=1$, we then obtain: 
\begin{eqnarray}
\label{eq:qp-approx}
\begin{array}{l}
p_\tau \approx \delta^{-1} \,,\\[.1in]
q_\tau \approx \delta^{2{\rm ht}(\alpha)-1} E(\alpha) \,,
\end{array}
\end{eqnarray}
where, as before,
$\alpha$ is the positive root corresponding to~$\tau$. 
Now let the operators $\widetilde M_\tau$ be given by (\ref{eq:Mtilde})
with $\gamma_w=\delta^{\l(w)}$. 
Then
\[
\widetilde M_\tau(w) = \delta^{\l(\tau w)-\l(w)} M_\tau(w) \,.
\]
Combining this with (\ref{eq:qp-approx}) and (\ref{eq:M}) yields 
\begin{eqnarray*}
\begin{array}{rcl}
\label{eq:Mtilde-specialized}
\widetilde M_\tau(w) &\approx& \left\{ 
\begin{array}{ll}
\delta^{\l(\tau w)-\l(w)-1} \,\tau w  
   & \mbox{if\ $\l(\tau w)>\l(w)$\ ;} \\[.1in]
\delta^{\l(\tau w)-\l(w)+2{\rm ht}(\alpha)-1} E(\alpha)\, \tau w   
   & \mbox{if\ $\l(\tau w)<\l(w)$\ .} \\
\end{array} \right.
\end{array}
\end{eqnarray*}
Note that always $\l(\tau)\leq 2{\rm ht}(\alpha)-1$; hence 
both exponents of~$\delta$ in the last formula
are \emph{nonnegative}. 
Letting $\delta\to 0$, we obtain the \emph{quantum Bruhat operators} 
$
Q_\tau=\displaystyle\lim_{\delta\to 0}\widetilde M_\tau 
$
given by 
\begin{equation}
\label{eq:Q}
Q_\tau(w)
=  
\left\{ 
\begin{array}{ll}
\tau w  
   & \mbox{if\ $\l(\tau w)=\l(w)+1$\ ;} \\[.1in]
E(\alpha)\, \tau w   
   & \mbox{if\ $\l(\tau w)=\l(w)-\l(\tau)$\ 
and\ $\l(\tau)=2{\rm ht}(\alpha)-1$;} \\[.1in]
0  & \mbox{otherwise\ .}
\end{array} \right.
\end{equation}
For the symmetric group, the requirement $\l(\tau)=2{\rm
  ht}(\alpha)-1$
in (\ref{eq:Q}) is superfluous, and we recover the type~$A$ quantum
 Bruhat operators  of~\cite{FK}.

Since the operators $Q_\tau$ were obtained from the mixed Bruhat operators 
of Theorem~\ref{th:main-W}  by specializing parameters, 
rescaling, and taking a limit, we have arrived at the following result.

\begin{corollary}
\label{cor:Q-ybe}
Let $\{Q_\tau\}_{\tau\in T}$ be the quantum Bruhat operators
defined by~{\rm(\ref{eq:Q})}.
Then the operators~$R_\tau=1+\varepsilon Q_\tau$ 
satisfy the Yang-Baxter equations
{\rm (\ref{eq:ybe-a1})--(\ref{eq:ybe-a2})}. 
\end{corollary}

For the type~$A$, it was noted in~\cite{FK} 
that the operators~$Q_\tau$ satisfy
the \emph{classical} Yang-Baxter equation~(\ref{eq:cybe}), 
which is a slightly weaker statement than Corollary~\ref{cor:Q-ybe}.

\medskip

We will now briefly explain the connection between our quantum Bruhat
operators and  the quantum cohomology of the generalized flag
manifold~$G/B$. 
Here $G$ is a semisimple connected complex Lie group associated with
the \emph{dual} root system~$\Phi^{\vee}$,
and $B$ is a Borel subgroup in~$G$. 
Let us identify each element $w\in W$ with the 
\emph{Schubert class}
\[
[w]=\sigma_{w^{-1}}
=[\,\overline{(Bw^{-1}B)/B}\,]
\in {\rm H}^{2\l(w)}(G/B,\mathbb{Z}) \,,
\]
 viewed 
as an element of the small quantum cohomology ring. 
(The reader is referred to \cite{fulton,FGP} and references therein
for relevant background.)  
In particular, the generators $s\in S$ will correspond to \emph{special}
Schubert classes~$[s]$.  
Extending the map $w\mapsto [w]$ to a linear isomorphism 
between $\k[W]$ and the (quantum) cohomology ring
assigns obvious meaning to expressions of the form $[Q(w)]$,
where $Q$ is an operator acting in~$\k[W]$. 

Let the quantum Bruhat operators $Q_\tau$ be given by~(\ref{eq:Q}), 
where the values of the multiplicative function $\alpha\mapsto
E(\alpha)$ at simple roots are set equal to the corresponding deformation
parameters of the quantum cohomology ring. 
Then the (quantum Monk's) formula for quantum multiplication of an
arbitrary Schubert class 
$[w]$ by a special Schubert class $[s]$ can be written as
follows: 
\begin{eqnarray}
\label{eq:q-monk}
\begin{array}{rcl}
 [w]*[s] &=& 
\quad\displaystyle\sum_{\alpha>0} \,\langle \omega, \alpha\rangle \,
[Q_\tau(w)] \\[.2in]
&=&   \displaystyle\mathop{\sum_{\alpha>0}}_{\l(w\tau)=\l(w)+1}
\langle \omega, \alpha\rangle\, [w\tau]
\,+ \displaystyle\mathop{\sum_{\alpha>0}}_{\l(w\tau)=\l(w)-2{\rm
    ht}(\alpha)+1} 
\langle \omega, \alpha\rangle\, E(\alpha)\, [w\tau] \,,
\end{array}
\end{eqnarray}
where, as before, the reflection~$\tau$ corresponds
to the positive root~$\alpha\in\Phi$,
and $\omega$ denotes the fundamental weight corresponding to~$s$. 

For the type~$A$ case, formula (\ref{eq:q-monk}) was first stated and
 proved in~\cite{FGP}.
For a general type, it was given by D.~Peterson (reproduced
in~\cite{carrell}, without proof).

\section{Motivation and proof}
\label{sec:proof}

In this section, we prove Theorem~\ref{th:main-W}, 
and also explain the origin of our solution 
(\ref{eq:global-p})--(\ref{eq:global-q}). 

Let us investigate the problem of choosing the parameters $p_\tau$ and
$q_\tau$ so that the operators $R_\tau$ given by
(\ref{eq:M})--(\ref{eq:R}) satisfy the Yang-Baxter
equations~(\ref{eq:ybe}). 
First of all, one easily checks that, for any choice of parameters, 
operators $R_\tau$  and $R_\sigma$
commute whenever $\tau$ and $\sigma$ do.
Therefore we only need to take care of~(\ref{eq:ybe}) in the cases
where both sides involve at least three factors.
In particular, for type~$A$ we only have to make sure 
that the operators $R_{ij}$ 
satisfy the quantum Yang-Baxter equation~(\ref{eq:ybe-a2}).

\subsection{Cosets modulo dihedral subgroups}
Notice that each operator $R_\tau$ stabilizes the span of every left coset
$W'w$ for any subgroup $W'$ containing~$\tau$. 
Let $W'$ be a dihedral reflection subgroup. 
(Thus $W'$ is of type $A_2\,$, $B_2\,$, or~$G_2\,$.) 
Then the span of every left coset of ~$W'$ 
is invariant under all operators appearing in the corresponding
Yang-Baxter equation~(\ref{eq:ybe}).
Thus the operators~$R_\tau$ satisfy~(\ref{eq:ybe}) if and
only if so do the restrictions of these operators onto
each space $\k[W'w]$ (which has dimension 6, 8, or~12).
Our plan is to explicitly write down the matrices of these
restrictions, plug them into the Yang-Baxter equation, 
and derive the complete set of equations for the parameters 
$\{p_\tau\}$ and  $\{q_\tau\}$.

The first step is to understand the combinatorics of the coset $W'w$ 
as a subposet of the Bruhat order. 
The following statement is known to hold for any Coxeter
group (see Dyer~\cite{Dy2}); in the special case of a Weyl group, 
it has a simple proof provided below. 

\begin{lemma}
\label{lem:cosets} 
Let $W'$ be a reflection subgroup of~$W$,
and let $S'$ be its set of canonical generators. 
Then the Bruhat order on $W'$ (viewed as a Coxeter group with
generating set~$S'$) 
coincides with the partial order induced from the Bruhat order on~$W$. 

With respect to the Bruhat order on $W$,
each coset $W'w$ has a unique minimal element~$\tilde w$.
For any $w'\in W'$ and $t\in T\cap W'$, 
we have $\l(t w')<\l(w')$
if and only if $\l(t w'\tilde w)<\l(w'\tilde w)$. 
\end{lemma}

\proof 
For $t\in T$, $w\in W$, the condition $\l(tw)<\l(w)$
is equivalent to $w^{-1}(\alpha)<0$, where $\alpha$ is the positive root
associated with~$t$. 
This implies the first part of the lemma.
To prove the second part, choose $\tilde w$ to be the element of
minimal length in~$W'w$ (if there are several such, pick any).
Take any reflection $t\in T\cap W'$ and the corresponding positive
root~$\alpha$. 
Then $\l(t\tilde w)>\l(\tilde w)$ and therefore
${\tilde w}^{-1}(\alpha)>0$. 
Thus ${\tilde w}^{-1}$ maps every positive root that
corresponds to a reflection in~$W'$ into a positive root
(and every negative into a negative).
Hence
\[
\l(t w')<\l(w') \Longleftrightarrow
(w')^{-1}(\alpha)<0 \Longleftrightarrow
{\tilde w}^{-1}(w')^{-1}(\alpha)<0 \Longleftrightarrow
\l(t w'\tilde w)<\l(w'\tilde w) \ ,
\]
as desired. 
\endproof

\begin{remark}
\label{remark:action}{\rm
Let $W'$ be a dihedral subgroup of~$W$.
The second part of Lemma~\ref{lem:cosets} implies 
that the action of the mixed
Bruhat operators participating in the Yang-Baxter equation for $W'$
restricted to each invariant subspace $\k[W'w]$ 
is canonically isomorphic to 
their action on $\k[W']$ 
via the linear isomorphism $w'\mapsto w'\tilde w$,
where $\tilde w$ is the unique minimal element of~$W'w$. 
In turn, the action on $\k[W']$ can be described quite explicitly
using the first part of Lemma~\ref{lem:cosets}:
the operators act as if $W'$ was the whole group. 
}
\end{remark}

\subsection{Example: solution for the symmetric group}
Let $W$ be the symmetric group $S_n\,$. 
For the convenience of the reader (and also to motivate subsequent
constructions), we will first treat this special case in complete
detail, and  later use it as a prototype for the general case. 

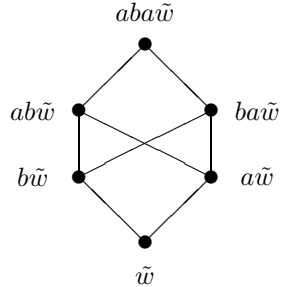
\begin{figure}[ht]
\setlength{\unitlength}{2.5pt} 
\begin{center}
\begin{picture}(40,40)(0,7)
\put(20,10){\line(1,1){10}}
\put(20,10){\line(-1,1){10}}
\put(10,20){\line(0,1){10}}
\put(10,30){\line(1,1){10}}
\put(20,40){\line(1,-1){10}}
\put(30,30){\line(0,-1){10}}
\put(10,20){\line(2,1){20}}
\put(10,30){\line(2,-1){20}}
\put(20,10){\circle*{2}}
\put(10,20){\circle*{2}}
\put(10,30){\circle*{2}}
\put(20,40){\circle*{2}}
\put(30,30){\circle*{2}}
\put(30,20){\circle*{2}}
\put(20,5){\makebox(0,0){$\tilde w$}}
\put(20,45){\makebox(0,0){$aba\tilde w$}}
\put(37,20){\makebox(0,0){$a\tilde w$}}
\put(37,30){\makebox(0,0){$ba\tilde w$}}
\put(3,20){\makebox(0,0){$b\tilde w$}}
\put(3,30){\makebox(0,0){$ab\tilde w$}}
\end{picture}
\end{center}
\caption{The Bruhat order on the coset $W'\tilde w$ in the symmetric
  group} 
\label{fig:coset}
\end{figure}

Let $W'$ be the 6-element 
dihedral reflection subgroup of $W=S_n$ generated by 
the reflections $a=(ij)$ and~$b=(jk)$, 
$1\leq i<j<k\leq n$.
A left coset of~$W'$ consists of the elements 
\begin{eqnarray}
\label{eq:6-dim}
\begin{array}{rcc}
\tilde w             &=& \cdots i \cdots j \cdots k \cdots \ ,\\[.05in]
a\tilde w            &=& \cdots j \cdots i \cdots k \cdots \ ,\\[.05in]
b\tilde w            &=& \cdots i \cdots k \cdots j \cdots \ ,\\[.05in]
ab\tilde w           &=& \cdots j \cdots k \cdots i \cdots \ ,\\[.05in]
ba\tilde w           &=& \cdots k \cdots i \cdots j \cdots \ ,\\[.05in]
aba\tilde w=bab\tilde w     &=& \cdots k \cdots j \cdots i \cdots \ ,\\
\end{array}
\end{eqnarray}
where all entries besides $i$, $j$, and $k$ are as in~$\tilde w$;
here $\tilde w$ is the minimal element of the coset. 
The partial order induced on $W'\tilde w$ from the Bruhat order on
$S_n$ is canonically isomorphic to the Bruhat order on the symmetric
group of permutations of three elements $i$, $j$, and~$k$. 
See Figure~\ref{fig:coset}. 

The restrictions of the operators $M_{ij}\,$, $M_{ik}\,$, and $M_{jk}$
to the invariant 6-dimen\-sional subspace $\k[W'\tilde w]$ spanned by the 
permutations~(\ref{eq:6-dim}) is readily computed using
Definition~\ref{def:mbo} and Remark~\ref{remark:action}. 
For example, let us compute $M_{ik}(ba\tilde w)$.
We have: $(ik)=bab$, 
$(ik)\cdot ba\tilde w =b\tilde w$,
$\l(b\tilde w)<\l(ba\tilde w)$,
implying $M_{ik}(ba\tilde w)=q_{ik}\,b\tilde w$. 
Analogous considerations show that in the linear basis of $\k[W'\tilde w]$ 
formed by the elements $\tilde w,a\tilde w,b\tilde w,ab\tilde
w,ba\tilde w,aba\tilde w$ (in this order), 
the restrictions of the operators 
$M_{ij}\,$, $M_{ik}\,$, and $M_{jk}$ 
are given by the following matrices:
\begin{eqnarray}
\label{eq:Mij-Mik}
\begin{array}{c}
M_{ij}=
\left[\!
\begin{array}{cccccc}
0 & \!q_{ij}\! & 0 & 0 & 0 & 0 \\
\!p_{ij}\! & 0 & 0 & 0 & 0 & 0 \\
0 & 0 & 0 & \!q_{ij}\! & 0 & 0 \\
0 & 0 & \!p_{ij}\! & 0 & 0 & 0 \\
0 & 0 & 0 & 0 & 0 & \!q_{ij}\! \\
0 & 0 & 0 & 0 & \!p_{ij}\! & 0 \\
\end{array}
\!\right],
\ 
M_{jk}=
\left[\!
\begin{array}{cccccc}
0 & 0 & \!q_{jk}\! & 0 & 0 & 0 \\
0 & 0 & 0 & 0 & \!q_{jk}\! & 0 \\
\!p_{jk}\! & 0 & 0 & 0 & 0 & 0 \\
0 & 0 & 0 & 0 & 0 &\! q_{jk}\! \\
0 & \!p_{jk}\! & 0 & 0 & 0 & 0 \\
0 & 0 & 0 & \!p_{jk}\! & 0 & 0 \\
\end{array}
\!\right],
\\[.6in]
M_{ik}=
\left[\begin{array}{cccccc}
0 & 0 & 0 & 0 & 0 & \!q_{ik}\! \\
0 & 0 & 0 & \!q_{ik}\! & 0 & 0 \\
0 & 0 & 0 & 0 &\! q_{ik}\! & 0\\
0 & \!p_{ik}\! & 0 & 0 & 0 & 0 \\
0 & 0 & \!p_{ik}\! & 0 & 0 & 0 \\
\!p_{ik}\! & 0 & 0 & 0 & 0 & 0 \\
\end{array}\right].
\end{array}
\end{eqnarray}
We are now prepared to write the conditions under which 
the operators $R_\tau = 1 + \varepsilon M_\tau$ 
satisfy the type~$A$ Yang-Baxter equation~(\ref{eq:ybe-a2}). 
The terms of degrees~0 and~1 in~$\varepsilon$ are clearly the same on
both sides of~(\ref{eq:ybe-a2}). Equating the quadratic terms gives 
the \emph{classical Yang-Baxter equation}~\cite{ybe}
\begin{equation}
\label{eq:cybe}
[M_{ij}\,,M_{jk}]=[M_{jk}\,,M_{ik}]+[M_{ik}\,,M_{ij}]
\end{equation}
(here $[A,B]=AB-BA$ stands for the commutator), 
while equating the cubic terms gives the quantum Yang-Baxter equation
for the~$M_\tau$: 
\begin{equation}
\label{eq:ybe-M}
M_{ij}M_{ik}M_{jk} = M_{jk}M_{ik}M_{ij}\quad 
\textrm{ if $\,i<j<k\,$.}
\end{equation}
Substituting (\ref{eq:Mij-Mik}) into (\ref{eq:cybe}), we obtain, upon
simplifications, the following system of equations: 
\begin{eqnarray}
\label{eq:cybe-6}
\left\{
\begin{array}{rcl}
- { q_{ij}}{ q_{jk}} + { p_{jk}}{ q_{ik}} + { q_{ik}}{ q_{ij}} &=&0\ ; \\[.1in]
  { q_{ij}}{ q_{jk}} - { q_{jk}}{ q_{ik}} - { q_{ik}}{ p_{ij}} &=&0\ ; \\[.1in]
  { p_{ij}}{ q_{jk}} - { q_{jk}}{ p_{ik}} - { q_{ik}}{ p_{ij}} &=&0\ ; \\[.1in]
- { q_{ij}}{ p_{jk}} + { p_{jk}}{ q_{ik}} + { p_{ik}}{ q_{ij}} &=&0\ ; \\[.1in]
  { p_{ij}}{ p_{jk}} - { q_{jk}}{ p_{ik}} - { p_{ik}}{ p_{ij}} &=&0\ ; \\[.1in]
- { p_{ij}}{ p_{jk}} + { p_{jk}}{ p_{ik}} + { p_{ik}}{ q_{ij}} &=&0\ . 
\end{array}
\right.
\end{eqnarray}
Making the same substitution into (\ref{eq:ybe-M}), we obtain
a single equation
$q_{ij}p_{ik}q_{jk}=p_{ij}q_{ik}p_{jk}\,$,
which actually follows from~(\ref{eq:cybe-6});
indeed, multiply the first equation in~(\ref{eq:cybe-6}) by~$p_{ik}\,$,
the last one---by~$q_{ik}\,$, and subtract. 

We thus arrived at the following result.

\begin{proposition}
\label{pro:system-6}
The operators $R_{ij}$ given by {\rm (\ref{eq:RA})}
satisfy the type~$A$ quantum Yang-Baxter equations~{\rm
  (\ref{eq:ybe-a1})--(\ref{eq:ybe-a2})} 
if and only if the parameters $\{p_{ij}\}$ and $\{q_{ij}\}$ 
satisfy the equations~{\rm (\ref{eq:cybe-6})}. 
\end{proposition}

It is possible to use equations~(\ref{eq:cybe-6})
to provide a complete parametric description of all solutions
of the Yang-Baxter equations of type~$A$ that have the form~(\ref{eq:R}).
However, this description is quite cumbersome
because of the many degenerate cases where lots of parameters
$p_{ij}$ and $q_{ij}$ vanish. 
Instead, we will now explicitly describe the particularly simple
family of solutions that is obtained in the ``generic'' case. 

Suppose for a moment that $q_{ij}\neq 0$ for any $i$ and~$j$.
Adding the first two equations in~(\ref{eq:cybe-6}) and dividing
by~$q_{ik}\,$, we obtain $p_{ij}-q_{ij}=p_{jk}-q_{jk}\,$.
Similarly, the second and third equations lead to
$p_{ij}-q_{ij}=p_{ik}-q_{ik}\,$.
This observation prompts the following consideration. 

Let us assume that the parameters $p_{ij}$ and $q_{ij}$ are related by 
\begin{equation}
\label{eq:p=q+c}
p_{ij} = q_{ij}+\kappa\ ,
\end{equation}
where $\kappa$ is a scalar constant that does not depend on $i$ and~$j$.
This assumption (motivated in the preceding paragraph)
immediately leads to substantial simplifications:
substituting (\ref{eq:p=q+c}) into (\ref{eq:cybe-6}) 
reduces this system of equations to a single equation
\begin{equation}
\label{eq:q(q+q+c)=qq}
q_{ik}(q_{ij}+q_{jk}+\kappa)= q_{ij} q_{jk} 
\end{equation}
---or, if you like, to
\begin{equation}
\label{eq:p(p+p-c)=pp}
p_{ik}(p_{ij}+p_{jk}-\kappa)= p_{ij} p_{jk} \,. 
\end{equation}

We conclude that whenever (\ref{eq:p=q+c}) and
(\ref{eq:q(q+q+c)=qq}) are satisfied by the collections of parameters
$p_{ij}$ and $q_{ij}\,$, the corresponding mixed Bruhat
operators of type~$A$ give rise to a solution of the Yang-Baxter
equations. 

Let us denote $q_i=q_{i,i+1}$ and $p_i=p_{i,i+1}=q_{i,i+1}+\kappa$ and
then use 
(\ref{eq:p=q+c})--(\ref{eq:p(p+p-c)=pp}) to compute all the $q_{ij}$
and $p_{ij}$ recursively. The prototypical example is $W=S_3\,$, in
which case 
we have 
\[
\begin{array}{l}
q_{13}=\displaystyle\frac{q_1 q_2}{q_1+q_2+\kappa}
=\frac{\kappa q_1 q_2}{p_1 p_2 - q_1 q_2} \,,\\[.25in]
p_{13}=\displaystyle\frac{p_1 p_2}{p_1+p_2-\kappa}
=\frac{\kappa p_1 p_2}{p_1 p_2 - q_1 q_2} \,.
\end{array}
\]
Continuing in the same fashion leads us to the formulas 
(\ref{eq:global-pA})--(\ref{eq:global-qA});
once those formulas are written down, proving them by induction on
$j-i$ is a matter of routine verification.

\begin{remark}
\label{rem:Sn-rep} 
{\rm
Observe that, for \emph{any} choice of parameters~$p_{ij}$
and~$q_{ij}$, the operators~$M_{ij}$ defined by~(\ref{eq:M})
are, up to a scalar, involutions: $M_{ij}^2=p_{ij}q_{ij}\,$. 
Now suppose that the $q_{ij}$ and $p_{ij}$ are given by
(\ref{eq:global-pA})--(\ref{eq:global-qA}). 
Then the~$M_{ij}$ (hence the normalized elements
$(p_{ij}q_{ij})^{-1/2}M_{ij}$) 
satisfy the Yang-Baxter equation 
$M_{ij}M_{ik}M_{jk} = M_{jk}M_{ik}M_{ij}$---but 
not the braid relation
\[
M_{i,i+1}M_{i+1,i+2}M_{i,i+1}=M_{i+1,i+2}M_{i,i+1}M_{i+1,i+2}\,.
\]
However, one can check that in the special case $q_1=\cdots=q_{n-1}$
the latter condition is satisfied, and we therefore obtain
a representation of the symmetric group~$S_n\,$. 
}\end{remark}

\subsection{Proof of Theorem~\ref{th:main-W}}

We will use our type $A$ solution as a model.
Let $W'$ be a dihedral subgroup of~$W$,
with the set of canonical generators $S'=\{a,b\}$. 
Thus $(W',S')$ is a Coxeter system of type~$A_2\,$, $B_2\,$,
or~$G_2\,$. 

Let $T'=T\cap W'$ be the set of reflections $\tau\in W'$;
these reflections correspond to the operators~$R_\tau$ involved in
the Yang-Baxter equation (\ref{eq:ybe}) associated with~$W'$. 
For any left coset $W'w$,  the subspace $\k[W'w]$ is invariant under the
action of all operators $R_\tau$ with $\tau\in T'$. 
By Lemma~\ref{lem:cosets}, the coset $W'w$ is in 
canonical bijection with~$W'$, giving rise to a canonically labelled
basis in the subspace~$\k[W'w]$.  
Furthermore, the matrices of the 
operators $R_\tau\,,\ \tau\in T'$, restricted 
to $\k[W'w]$, do not depend on the choice of a coset. 
We can explicitly write down these matrices (of size $6\times 6$,
$8\times 8$, or $12\times 12$) in terms of the corresponding parameters
$p_\tau$ and $q_\tau\,$, in complete analogy with (\ref{eq:Mij-Mik}). 

If $W'$ is of type~$A_2\,$, then we obtain the matrices whose only
difference from (\ref{eq:Mij-Mik}) is in notation:
we have to replace the subscripts $ij$, $jk$, and $ik$ 
by the reflections $a$, $b$ and $aba=bab$, respectively. 
This leads to a system of equations of the form (\ref{eq:cybe-6}).
In view of (\ref{eq:p=q+c}) and (\ref{eq:q(q+q+c)=qq}), 
these equations will be satisfied if we impose the condition 
\begin{equation}
\label{eq:general-p=q+c}
p_\tau = q_\tau + \kappa_\tau \ , \ \ \tau\in T\ ,
\end{equation}
where $\kappa_\tau$ only depends on whether $\tau$ corresponds to a long or
a short root, 
and require that 
\begin{equation}
\label{eq:general-q(q+q+c)=qq}
q_a q_b = q_{aba}(q_a+q_b+\kappa_a)
\end{equation}
whenever $a,b\in T$ are canonical generators for a dihedral
subgroup of type~$A_2\,$.

If $W'$ is of type $B_2\,$, then the set $T'$ consists of four reflections
$a$, $b$, $aba$, and $bab$.
Labelling the basis of the invariant 8-dimensional subspace 
by the elements of $W'=\{e,a,b,ba,ab,aba,bab,abab\}$ (in this order),
we obtain the matrices
\[
M_{a}=
\left[\!
\begin{array}{cccccccc}
0 &\!q_{a}\! & 0 & 0 & 0 & 0 & 0 & 0  \\
\!p_{a}\! & 0 & 0 & 0 & 0 & 0 & 0 & 0 \\
0 & 0 & 0 & 0 & \!q_{a}\! & 0 & 0 & 0 \\
0 & 0 & 0 & 0 & 0 & \!q_{a}\! & 0 & 0 \\
0 & 0 & \!p_{a}\! & 0 & 0 & 0 & 0 & 0 \\
0 & 0 & 0 & \!p_{a}\! & 0 & 0 & 0 & 0 \\
0 & 0 & 0 & 0 & 0 & 0 & 0 & \!q_{a}\! \\
0 & 0 & 0 & 0 & 0 & 0 & \!p_{a}\! & 0 \\
\end{array}
\!\right],
\ 
M_{b}=
\left[\!
\begin{array}{cccccccc}
0 & 0 & \!q_{b}\! & 0 & 0 & 0 & 0 & 0  \\
0 & 0 & 0 & \!q_{b}\! & 0 & 0 & 0 & 0 \\
\!p_{b}\! & 0 & 0 & 0 & 0 & 0 & 0 & 0 \\
0 & \!p_{b}\! & 0 & 0 & 0 & 0 & 0 & 0 \\
0 & 0 & 0 & 0 & 0 & 0 & \!q_{b}\! & 0 \\
0 & 0 & 0 & 0 & 0 & 0 & 0 & \!q_{b}\! \\
0 & 0 & 0 & 0 & \!p_{b}\! & 0 & 0 & 0 \\
0 & 0 & 0 & 0 & 0 & \!p_{b}\! & 0 & 0 \\
\end{array}
\!\right],
\]
and, in a similar way, the matrices $M_{aba}$ and $M_{bab}\,$.
Substituting these matrices into the type~$B_2$ Yang-Baxter equation
\[
(1+M_a)(1+M_{aba})(1+M_{bab})(1+M_b)=
(1+M_b)(1+M_{bab})(1+M_{aba})(1+M_a),
\]
we obtain a system of equations for the 8 parameters $p_\tau$ and
$q_\tau$ corresponding to $\tau\in\{a,b,aba,bab\}$.
If we make an assumption~(\ref{eq:general-p=q+c}),
this system of equations collapses into the single equation
\begin{equation}
\label{eq:cybe-B}
q_a q_b = q_a q_{aba} + q_{aba} q_{bab} + q_{bab} q_b
         + \kappa_a q_{aba} + \kappa_b q_{bab} \ ,
\end{equation}
which we want to be satisfied whenever $a,b\in T$ are canonical
generators for a dihedral subgroup of type~$B_2\,$. 

For type~$G_2\,$, we have 12 parameters $p_\tau$ and $q_\tau\,$.
Assuming~(\ref{eq:general-p=q+c}), we express everything in terms of
the 6 parameters $q_\tau$ and the 2 parameters $\kappa_\tau$ (for the short
and long roots, respectively).
Making a substitution into the Yang-Baxter equation of type~$G_2\,$,
we obtain the two equations 
\begin{eqnarray}
\label{eq:cybe-G1}
\begin{array}{l}
q_a q_b =  q_a q_{aba} + q_{aba} q_{ababa} + q_{ababa} q_{babab} +
         q_{babab} q_{bab} + q_{bab} q_b \\[.1in]
\qquad \qquad          
+ \kappa_a q_{aba} + \kappa_a q_{babab} + \kappa_b q_{bab} + \kappa_b
         q_{ababa} 
\end{array}
\end{eqnarray}
and
\begin{eqnarray}
\label{eq:cybe-G2}
\begin{array}{l}
-q_a q_{bab}
+q_a q_{ababa}
-q_b q_{aba}
+q_b q_{babab} 
+q_{aba} q_{babab} 
+q_{bab} q_{ababa} \\[.1in]

+q_a q_b q_{aba} q_{bab} 
-q_a q_b q_{aba} q_{ababa}
-q_a q_b q_{bab} q_{babab}
-q_a q_b q_{ababa} q_{babab} \\[.1in]

+q_a q_{aba} q_{bab} q_{babab}
+q_b q_{bab} q_{ababa} q_{babab}
+q_b q_{aba} q_{bab} q_{ababa} \\[.1in]
+q_a q_{aba} q_{ababa} q_{babab} 
+q_{aba} q_{bab} q_{ababa} q_{babab}  \\[.1in]

+\kappa_a (q_{ababa}
-q_a q_b q_{babab} 
+q_b q_{aba} q_{bab} 
+q_a q_{aba} q_{babab}  \\[.1in]
\qquad +q_{aba} q_{bab} q_{babab} 
+q_{aba} q_{ababa} q_{babab}) \\[.1in]

+\kappa_b (q_{babab}
-q_a q_b q_{ababa} 
+q_a q_{aba} q_{bab} 
+q_b q_{bab} q_{ababa}  \\[.1in]
\qquad +q_{aba} q_{bab} q_{ababa} 
+q_{bab} q_{ababa} q_{babab}) \\[.1in]

+\kappa_a^2 q_{aba} q_{babab} 
+\kappa_a \kappa_b q_{aba} q_{bab} 
+\kappa_b^2 q_{bab} q_{ababa}

= 0 \,. \\[.1in]
\end{array}
\end{eqnarray}

We are now fully prepared to complete the proof of
Theorem~\ref{th:main-W},
which amounts to checking the equations
(\ref{eq:general-q(q+q+c)=qq}), (\ref{eq:cybe-B}), and
(\ref{eq:cybe-G1})--(\ref{eq:cybe-G2})
for every dihedral subgroup $W'$ of type $A_2\,$, $B_2\,$, or
$G_2\,$, respectively,
provided the $q_\tau$ and $p_\tau$ are given by
(\ref{eq:global-p})--(\ref{eq:global-q}).  
This is a straightforward verification.
Let $a$ and $b$ be the canonical generators of~$W'$,
and let $\alpha$ and $\beta$ be the corresponding positive roots. 
For $W'$ of type~$A_2\,$, we have $\kappa_a=\kappa_b\,$, and equation
(\ref{eq:general-q(q+q+c)=qq}) becomes
\[
\begin{array}{l}
\frac{E_2(\alpha)}{E_1(\alpha)-E_2(\alpha)}\,\,
\frac{E_2(\beta)}{E_1(\beta)-E_2(\beta)}\\[.2in]
=
\frac{E_2(\alpha+\beta)}
{E_1(\alpha+\beta)-E_2(\alpha+\beta)}\,
\left(
\frac{E_2(\alpha)}{E_1(\alpha)-E_2(\alpha)}\,+\,
\frac{E_2(\beta)}{E_1(\beta)-E_2(\beta)}
+1
\right)\, ,
\end{array}
\]
which is immediately checked using that $E_1$ and $E_2$ are
multiplicative (cf.~(\ref{eq:multiplicativity})). 
Let $W'$ be of type~$B_2\,$.
Note that equation  (\ref{eq:cybe-B}) 
is invariant under interchanging $a$ and~$b$. 
Therefore without loss of generality we may assume that 
$\alpha$ is short while $\beta$ is long.
Then $aba$ and $bab$ correspond to positive roots
$2\alpha+\beta$ (long) and $\alpha+\beta$ (short), respectively. 
Substituting (\ref{eq:global-q}) into (\ref{eq:cybe-B}) and factoring
out $\kappa_a \kappa_b\,$, we obtain 
\[
\begin{array}{l}
\frac{E_2(\alpha)}{E_1(\alpha)-E_2(\alpha)}\,
\frac{E_2(\beta)}{E_1(\beta)-E_2(\beta)}\\[.2in]
=
\frac{E_2(\alpha)}{E_1(\alpha)-E_2(\alpha)}\,
\frac{E_2(2\alpha+\beta)}
{E_1(2\alpha+\beta)-E_2(2\alpha+\beta)}
\,+\,
\frac{E_2(2\alpha+\beta)}
{E_1(2\alpha+\beta)-E_2(2\alpha+\beta)}\,
\frac{E_2(\alpha+\beta)}
{E_1(\alpha+\beta)-E_2(\alpha+\beta)} \\[.2in]
\,+\,
\frac{E_2(\alpha+\beta)}
{E_1(\alpha+\beta)-E_2(\alpha+\beta)}\,
\frac{E_2(\beta)}{E_1(\beta)-E_2(\beta)}
\,+\,
\frac{E_2(2\alpha+\beta)}
{E_1(2\alpha+\beta)-E_2(2\alpha+\beta)}\,+\,
\frac{E_2(\alpha+\beta)}
{E_1(\alpha+\beta)-E_2(\alpha+\beta)}\,,
\end{array}
\]
which is easily checked using~(\ref{eq:multiplicativity}). 
The case $G_2$ is verified in a similar way
(preferably with the help of a computer). 
\endproof

\section{
Tilted Bruhat orders 
}
\label{sec:combin} 

We will now apply the results of Section~\ref{sec:Q}
to the combinatorics of the Coxeter system~$(W,S)$. 
Our main tool will be the following specialization of quantum Bruhat
operators~(\ref{eq:Q}). 

\begin{corollary}
\label{cor:Qq}
Let 
\begin{equation}
\label{eq:Qq}
Q_\tau(w)
=  
\left\{ 
\begin{array}{ll}
\tau w  
   & \mbox{if\ $\l(\tau w)=\l(w)+1$\ ;} \\[.1in]
\tau w   
   & \mbox{if $\l(\tau w)=\l(w)-\l(\tau)$\ and\ 
$\l(\tau)=2{\rm ht}(\alpha)-1$;} \\[.1in]
0  & \mbox{otherwise\ .}
\end{array} \right.
\end{equation}
Then the operators 
\begin{equation}
\label{eq:Rq}
R_\tau = 1+\varepsilon Q_\tau
\end{equation}
satisfy the Yang-Baxter equations.
\end{corollary}

\proof
In Corollary~\ref{cor:Q-ybe}, set $E(\alpha)=1$.
\endproof

\begin{definition}
\label{def:D(W)}
{\rm 
Motivated by 
(\ref{eq:Qq}), 
let us introduce the following digraph~$D(W)$.  
The vertices of $D(W)$ are the elements of the group~$W$.
For $u\in W$ and $\tau\in T$, we put a directed edge
from $u$ to $v=\tau u$ if either $\l(v)=\l(u)+1$ 
or $\l(v)=\l(u)-\l(\tau)$ and $\l(\tau)=2{\rm ht}(\alpha)-1$,
where $\alpha$ is the corresponding positive root. 
In other words, $(u,\tau u)$ is an edge in $D(W)$
if multiplying $u$ by $\tau$ on the left 
either increases the length of $u$ by as little as possible
or decreases the length of $u$ by as much as possible,
given the height of~$\alpha$. 

Once a reflection ordering~$\varphi$ for~$W$ is chosen,
we label the edges of $D(W)$ by assigning label~$\varphi(\tau)$
to an edge~$(u,v)$ with $v=\tau u$.  
We will write $u\stackrel{m}\to v$ to denote that $(u,v)$ is
an edge in $D(W)$ labelled by~$m$.
}\end{definition}

\begin{figure}[ht]
\setlength{\unitlength}{1.8pt} 
\begin{center}
\begin{picture}(40,90)(10,-5)

\put(20.5,0.5){\vector(-1,1){10}}
\put(10.5,10.5){\line(-1,1){10}}
\put(-0.5,19.5){\vector(1,-1){10}}
\put(9.5,9.5){\line(1,-1){10}}

\put(19.5,0.5){\vector(1,1){10}}
\put(29.5,10.5){\line(1,1){10}}
\put(40.5,19.5){\vector(-1,-1){10}}
\put(30.5,9.5){\line(-1,-1){10}}

\put(39.5,59.5){\vector(-1,1){10}}
\put(29.5,69.5){\line(-1,1){10}}
\put(20.5,80.5){\vector(1,-1){10}}
\put(30.5,70.5){\line(1,-1){10}}

\put(0.5,59.5){\vector(1,1){10}}
\put(10.5,69.5){\line(1,1){10}}
\put(19.5,80.5){\vector(-1,-1){10}}
\put(9.5,70.5){\line(-1,-1){10}}

\put(-0.7,20){\vector(0,1){10}}
\put(-0.7,30){\line(0,1){10}}
\put(0.7,40){\vector(0,-1){10}}
\put(0.7,30){\line(0,-1){10}}

\put(-0.7,40){\vector(0,1){10}}
\put(-0.7,50){\line(0,1){10}}
\put(0.7,60){\vector(0,-1){10}}
\put(0.7,50){\line(0,-1){10}}

\put(39.3,20){\vector(0,1){10}}
\put(39.3,30){\line(0,1){10}}
\put(40.7,40){\vector(0,-1){10}}
\put(40.7,30){\line(0,-1){10}}

\put(39.3,40){\vector(0,1){10}}
\put(39.3,50){\line(0,1){10}}
\put(40.7,60){\vector(0,-1){10}}
\put(40.7,50){\line(0,-1){10}}

\put(0,20){\vector(2,1){28}}
\put(28,34){\line(2,1){12}}

\put(0,40){\vector(2,1){15}}
\put(15,47.5){\line(2,1){25}}

\put(40,20){\vector(-2,1){28}}
\put(12,34){\line(-2,1){12}}

\put(40,40){\vector(-2,1){15}}
\put(25,47.5){\line(-2,1){25}}

\put(40,60){\vector(-1,-3){10}}
\put(30,30){\line(-1,-3){10}}

\put(20,80){\vector(-1,-3){10}}
\put(10,50){\line(-1,-3){10}}

\put(20,0){\circle*{2}}
\put(0,20){\circle*{2}}
\put(0,40){\circle*{2}}
\put(0,60){\circle*{2}}
\put(20,80){\circle*{2}}
\put(40,20){\circle*{2}}
\put(40,40){\circle*{2}}
\put(40,60){\circle*{2}}

\put(20,-5){\makebox(0,0){$e$}}
\put(-7,20){\makebox(0,0){$a$}}
\put(-7,40){\makebox(0,0){$ba$}}
\put(-7,60){\makebox(0,0){$aba$}}
\put(47,20){\makebox(0,0){$b$}}
\put(47,40){\makebox(0,0){$ab$}}
\put(47,60){\makebox(0,0){$bab$}}
\put(20,85){\makebox(0,0){$\wnot$}}

\put(7,7){\makebox(0,0){\textbf{\small 1}}}
\put(33,7){\makebox(0,0){\textbf{\small 4}}}
\put(-4,30){\makebox(0,0){\textbf{\small 4}}}
\put(44,30){\makebox(0,0){\textbf{\small 1}}}
\put(-4,50){\makebox(0,0){\textbf{\small 1}}}
\put(44,50){\makebox(0,0){\textbf{\small 4}}}
\put(7,73){\makebox(0,0){\textbf{\small 4}}}
\put(33,73){\makebox(0,0){\textbf{\small 1}}}

\put(15,36.5){\makebox(0,0){\textbf{\small 3}}}
\put(25,56){\makebox(0,0){\textbf{\small 2}}}

\put(15,24){\makebox(0,0){\textbf{\small 2}}}
\put(25,44){\makebox(0,0){\textbf{\small 3}}}

\put(23,17){\makebox(0,0){\textbf{\small 3}}}
\put(17,63){\makebox(0,0){\textbf{\small 3}}}

\put(80,50){\mbox{\textbf{\small 1} \quad $a$}} 
\put(80,40){\mbox{\textbf{\small 2} \quad $aba$}} 
\put(80,30){\mbox{\textbf{\small 3} \quad $bab$}}
\put(80,20){\mbox{\textbf{\small 4} \quad $b$}}

\end{picture}
\end{center}
\caption{The digraph  $D(W)$ for $W$ of type $B_2$}
\label{fig:D(B_2)}
\end{figure}
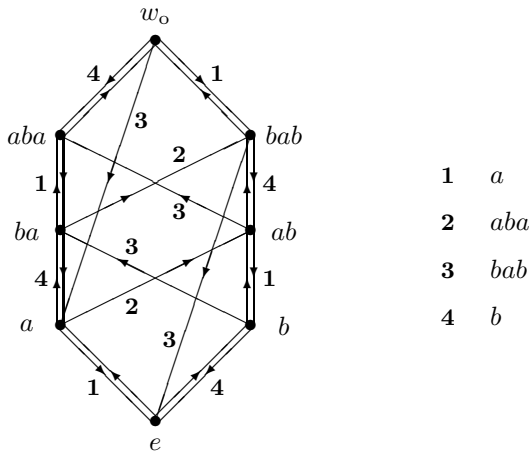

\begin{example}{\rm
Consider a Weyl group of type~$B_2\,$.
This is the first instance where the condition 
$\l(\tau)=2{\rm ht}(\alpha)-1$ comes into play.
Let $a$ and $b$ be the generators of~$W$ that correspond to the simple
roots $\alpha$ (short) and $\beta$ (long), respectively.
Then the reflections $a$, $b$, and $bab$ satisfy this condition,
while $aba$ does not (see Figure~\ref{fig:roots-B2}).

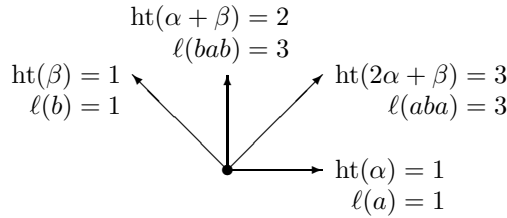
\begin{figure}[ht]
\setlength{\unitlength}{1.8pt} 
\begin{center}
\begin{picture}(40,40)(-10,-5)

\put(0,0){\circle*{2}}

\put(0,0){\vector(1,0){20}}
\put(0,0){\vector(1,1){20}}
\put(0,0){\vector(0,1){20}}
\put(0,0){\vector(-1,1){20}}

\put(18,-5){\mbox{
$\begin{array}{r}
{\rm ht}(\alpha)=1\\ 
\l(a)=1
\end{array}
$}}

\put(18,15){\mbox{
$\begin{array}{r}
{\rm ht}(2\alpha+\beta)=3\\ 
\l(aba)=3
\end{array}
$}}

\put(-25,27){\mbox{
$\begin{array}{r}
{\rm ht}(\alpha+\beta)=2\\ 
\l(bab)=3
\end{array}
$}}

\put(-50,15){\mbox{
$\begin{array}{r}
{\rm ht}(\beta)=1\\ 
\l(b)=1
\end{array}
$}}

\end{picture}
\end{center}
\caption{Root system $B_2$}
\label{fig:roots-B2}
\end{figure}

\noindent
We hence disallow down-directed edges that correspond
to multiplying by~$aba$ (on the left). 
The resulting graph $D(W)$, for the reflection ordering $a<aba<bab<b$,
is shown in Figure~\ref{fig:D(B_2)}. 
}
\end{example}

Notice that the construction of the digraph~$D(W)$ depends on the root
system~$\Phi$, not just on the Weyl group~$W$.
Thus, for example, digraphs of types $B_n$ and $C_n$ will differ from each
other. 

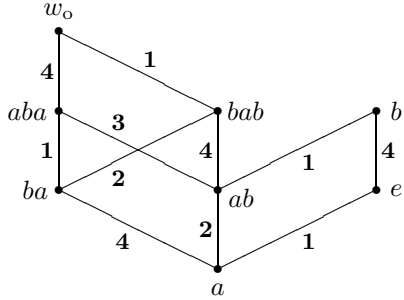
\begin{figure}[ht]
\setlength{\unitlength}{1.5pt} 
\begin{center}
\begin{picture}(40,70)(10,-5)

\put(40,0){\line(-2,1){40}}
\put(40,20){\line(-2,1){40}}
\put(40,40){\line(-2,1){40}}

\put(40,0){\line(2,1){40}}
\put(0,20){\line(2,1){40}}
\put(40,20){\line(2,1){40}}

\put(40,0){\line(0,1){20}}
\put(0,20){\line(0,1){20}}
\put(0,40){\line(0,1){20}}
\put(80,20){\line(0,1){20}}
\put(40,20){\line(0,1){20}}

\put(40,0){\circle*{2}}
\put(0,20){\circle*{2}}
\put(40,20){\circle*{2}}
\put(80,20){\circle*{2}}
\put(0,40){\circle*{2}}
\put(40,40){\circle*{2}}
\put(80,40){\circle*{2}}
\put(0,60){\circle*{2}}

\put(40,-5){\makebox(0,0){$a$}}
\put(-6,20){\makebox(0,0){$ba$}}
\put(46,18){\makebox(0,0){$ab$}}
\put(85,20){\makebox(0,0){$e$}}
\put(-8,40){\makebox(0,0){$aba$}}
\put(47,40){\makebox(0,0){$bab$}}
\put(85,40){\makebox(0,0){$b$}}
\put(0,65){\makebox(0,0){$\wnot$}}

\put(63,7){\makebox(0,0){\textbf{\small 1}}}
\put(37,10){\makebox(0,0){\textbf{\small 2}}}
\put(16,7){\makebox(0,0){\textbf{\small 4}}}

\put(-3,30){\makebox(0,0){\textbf{\small 1}}}
\put(15,23){\makebox(0,0){\textbf{\small 2}}}
\put(15,37){\makebox(0,0){\textbf{\small 3}}}
\put(37,30){\makebox(0,0){\textbf{\small 4}}}

\put(83,30){\makebox(0,0){\textbf{\small 4}}}
\put(63,27){\makebox(0,0){\textbf{\small 1}}}

\put(-3,50){\makebox(0,0){\textbf{\small 4}}}

\put(23,53){\makebox(0,0){\textbf{\small 1}}}

\end{picture}
\end{center}
\caption{Tilted Bruhat order $D_a(W)$ for $W$ of type~$B_2$}
\label{fig:D_a(W)}
\end{figure}

\begin{definition}
\label{def:tilted}
{\rm
For $u,v\in W$, let $\l(u,v)$ denote the length of the shortest path
in $D(W)$ from $u$ to~$v$.
In particular, $\l(e,v)=\l(v)$ is the usual length function,
where $e$ denotes the identity element in~$W$;
moreover, $\l(u,v)=\l(v)-\l(u)$ whenever $u\leq v$ in the Bruhat order.
The \emph{tilted Bruhat interval} between $u$
and~$v$ is the set
\[
D(u,v) = \{ w\in W \,:\, \l(u,w)+\l(w,v)=\l(u,v) \} \,,
\]
equipped with the following partial order: $w_1\preceq w_2$ if and only if
\[
\l(u,w_1)+\l(w_1,w_2)+\l(w_2,v)=\l(u,v) \,.
\]
Thus $D(u,v)$ is a graded poset whose Hasse diagram is the
minimal subgraph of $D(W)$ containing all directed paths from $u$ to
$v$ that have the smallest possible length. 
If $u\leq v$ in the Bruhat order, then $D(u,v)$ is nothing but the
interval $[u,v]=\{w\,:\,u\leq w\leq v\}$, explaining our choice of
terminology. 
Note that the intervals $D(u,v)$ and $D(v,u)$ are by no means dual
posets; for example, in Figure~\ref{fig:D(B_2)} the interval
$D(e,\wnot)$ is the whole Bruhat order, while $D(\wnot,e)$ has only 
four vertices (see Figure~\ref{fig:D(ab,a)}). 

Let us also define the \emph{tilted Bruhat order} $D_u(W)$ as a
graded partial order on~$W$ and the following order relation:
$w_1\preceq_u w_2$ if and only if
\[
\l(u,w_1)+\l(w_1,w_2)=\l(u,w_2) \,.
\]
Thus $w_1\preceq_u w_2$ if and only if there exists a shortest path from
$u$ to $w_2$ 
that passes through~$w_1\,$. 
Note that any interval in this poset (or in any $D(u,v)$) is again a
tilted Bruhat interval between corresponding vertices. 
$D_e(W)$ is the usual Bruhat order. 

Any choice of reflection ordering induces edge labelling of the Hasse
diagrams of $D(u,v)$ and $D_u(W)$ inherited from~$D(W)$. 
Figure~\ref{fig:D(ab,a)} shows the tilted Bruhat
interval $D(ab,a)$ for $W$ of type~$B_2\,$, with the same conventions
as in Figure~\ref{fig:D(B_2)}. 
(It also happens to coincide with the tilted Bruhat
order~$D_{ab}(W)$.) 
Figure~\ref{fig:D_a(W)} shows an example of a tilted Bruhat
order which is not pure (i.e., does not have a~$\hat 1$). 
}\end{definition}

\begin{figure}[ht]
\setlength{\unitlength}{1.5pt} 
\begin{center}
\begin{picture}(140,70)(10,-5)

\put(40,0){\line(-2,1){40}}
\put(40,20){\line(-2,1){40}}
\put(80,20){\line(-2,1){40}}
\put(80,40){\line(-2,1){40}}

\put(40,0){\line(2,1){40}}
\put(0,20){\line(2,1){40}}
\put(40,20){\line(2,1){40}}
\put(0,40){\line(2,1){40}}

\put(40,0){\line(0,1){20}}
\put(0,20){\line(0,1){20}}
\put(80,20){\line(0,1){20}}
\put(40,40){\line(0,1){20}}

\put(40,0){\circle*{2}}
\put(0,20){\circle*{2}}
\put(40,20){\circle*{2}}
\put(80,20){\circle*{2}}
\put(0,40){\circle*{2}}
\put(40,40){\circle*{2}}
\put(80,40){\circle*{2}}
\put(40,60){\circle*{2}}

\put(40,-5){\makebox(0,0){$ab$}}
\put(-5,20){\makebox(0,0){$b$}}
\put(47,18){\makebox(0,0){$aba$}}
\put(88,20){\makebox(0,0){$bab$}}
\put(-7,40){\makebox(0,0){$ba$}}
\put(44,42){\makebox(0,0){$e$}}
\put(87,40){\makebox(0,0){$\wnot$}}
\put(40,65){\makebox(0,0){$a$}}

\put(63,8){\makebox(0,0){\textbf{\small 4}}}
\put(37,10){\makebox(0,0){\textbf{\small 3}}}
\put(16,8){\makebox(0,0){\textbf{\small 1}}}

\put(-3,30){\makebox(0,0){\textbf{\small 3}}}
\put(15,23){\makebox(0,0){\textbf{\small 4}}}
\put(15,37){\makebox(0,0){\textbf{\small 1}}}
\put(65,23){\makebox(0,0){\textbf{\small 3}}}
\put(65,37){\makebox(0,0){\textbf{\small 4}}}
\put(83,30){\makebox(0,0){\textbf{\small 1}}}

\put(63,53){\makebox(0,0){\textbf{\small 3}}}
\put(37,50){\makebox(0,0){\textbf{\small 1}}}
\put(16,53){\makebox(0,0){\textbf{\small 4}}}

\put(140,10){\line(-1,1){20}}
\put(140,10){\line(1,1){20}}
\put(120,30){\line(1,1){20}}
\put(160,30){\line(-1,1){20}}

\put(140,10){\circle*{2}}
\put(120,30){\circle*{2}}
\put(160,30){\circle*{2}}
\put(140,50){\circle*{2}}

\put(140,5){\makebox(0,0){$\wnot$}}
\put(115,30){\makebox(0,0){$a$}}
\put(167,30){\makebox(0,0){$bab$}}
\put(140,55){\makebox(0,0){$e$}}

\put(127,17){\makebox(0,0){\textbf{\small 3}}}
\put(153,17){\makebox(0,0){\textbf{\small 1}}}
\put(127,43){\makebox(0,0){\textbf{\small 1}}}
\put(153,43){\makebox(0,0){\textbf{\small 3}}}

\end{picture}
\end{center}
\caption{Tilted Bruhat intervals $D(ab,a)$ and $D(\wnot,e)$}
\label{fig:D(ab,a)}
\end{figure}
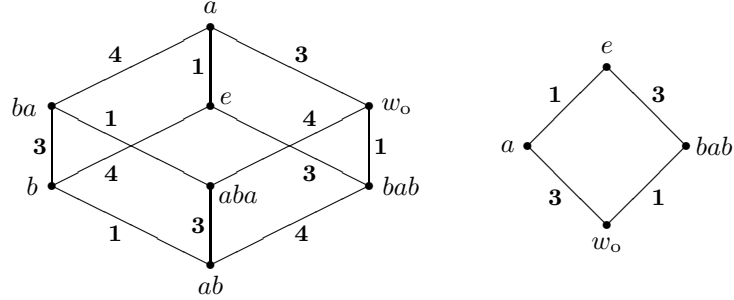

Our main combinatorial result is an extension of certain
fundamental property of Bruhat orders to their ``tilted analogues'' 
introduced in Definition~\ref{def:tilted}.
Let us first review the known facts.

Recall ~\cite{StaNato} that a finite graded poset with $\hat 0$ and
 $\hat 1$ (resp.\ with~$\hat 0$) 
 is called \emph{Eulerian} (resp.\ \emph{lower Eulerian}) 
if its \emph{M\"obius function}~\cite{rota} is given by
\[
\mu(x,y)=(-1)^{{\rm rank}(y)-{\rm rank}(x)}
\] for any $x\leq y$. 
A well known (but non-trivial---cf.~\cite{BW, DeoBru, KL})
theorem of Verma~\cite{V,V-correction} asserts that any interval in
the Bruhat order of any Coxeter group is Eulerian. 
To our knowledge, no simple proof of this result is known,
except for the special case $x=e$ (see
Lascoux~\cite[Lemma~1.13]{lascoux-grothendieck}). 
The story of Verma's theorem is described in \cite[p.~176]{Hum}. 
Remarkably, it can be strengthened as follows: 
any Bruhat interval is actually a face poset of a shellable regular CW
sphere 
(see~Bj\"orner~\cite[Theorem~5.1]{bjorner-cw} and
Bj\"orner-Wachs~\cite[Theorem~4.2]{BW}); 
hence it is also Cohen-Macaulay~\cite{BGS}.

All the statements mentioned in the preceding paragraph are implied by
the following  
``lexicographic shellability'' result conjectured by Bj\"orner and
proved by Dyer~\cite[Proposition~4.3]{Dy} for an arbitrary Coxeter
group. 
(This requires a more general definition of a reflection ordering,
not needed in  this paper.) 

\begin{theorem}
\label{EL}
\cite{Dy} Let $u,v \in W$, $u \leq v$. 
Then, for any reflection ordering, there exists
a unique label-increasing (and, by reversal of the ordering, also
unique label-decreasing) maximal chain from $u$ to~$v$ in the Bruhat
order of~$W$.
The sequence of labels associated with this chain is lexicographically
minimal (resp. lexicographically maximal) among all maximal chains
from $u$ to~$v$. 
\end{theorem}

We generalize this result (in the case of a Weyl group) as follows. 

\begin{theorem}
\label{th:tilted}
Fix a reflection ordering~$\varphi$ in a Weyl group~$W$. 
\nopagebreak
\medskip

\noindent
\textbf{1.}~For any pair of
elements $u,v\in W$, 
there is a unique path from $u$ to $v$ 
in the directed graph~$D(W)$ such that its sequence of labels is strictly
increasing (resp.\ strictly decreasing). 
\medskip

\noindent
\textbf{2.}~The unique label-increasing (resp.\ label-decreasing) path
from $u$ 
to~$v$ has the smallest possible length~$\l(u,v)$. 
Moreover it is lexicographically minimal (resp. lexicographically
maximal) among all shortest paths from $u$ to~$v$. 
\medskip

\noindent
\textbf{3.}~For any $u\in W$, 
\begin{equation}
\label{eq:Tu=sum-v}
R_{\varphi ^{-1}(1)}\cdots R_{\varphi ^{-1}(N)}(u) 
= \sum_{v \in W} \varepsilon^{\l(u,v)}\, v \,,  
\end{equation}
where the $R_\tau$ are given by 
{\rm (\ref{eq:Qq})--(\ref{eq:Rq});} as before, $N=\l(\wnot)$.  
\end{theorem}

\proof
We first note that part~1 of the theorem is equivalent to the special
case $\varepsilon=1$ of~(\ref{eq:Tu=sum-v}). 
Indeed, comparing our definition of the digraph~$D(W)$ to (\ref{eq:Qq}),
we see that $(u,v)$ is an edge in $D(W)$
if and only if $v=Q_\tau(u)$ for some $\tau \in T$,
in which case $(u,v)$ is labelled by~$\varphi(\tau)$.  
Thus the identity (\ref{eq:Tu=sum-v}), with $\varepsilon=1$, asserts 
existence and uniqueness of the label-decreasing path. 

Let us denote by $\T$ the specialization of the operator 
$R_{\varphi ^{-1}(1)}\cdots R_{\varphi ^{-1}(N)}$ obtained by setting
$\varepsilon=1$. 
By Proposition~\ref{pro:invariant-product} and Corollary~\ref{cor:Qq},
the operator $\T$
does not depend on the choice of reflection ordering~$\varphi$.

We will identify an element $w \in W$ 
with the linear operator $ u \mapsto wu$ in $\k[W]$.
Let $s\in S$.
Then (\ref{eq:Qq}) gives $Q_s=s$, implying
\begin{equation}
\label{eq:Rs=R}
(1+Q_s) s = 1+Q_s \ . 
\end{equation}
Since there exists a reduced decomposition of $\wnot$ that ends
in~$s$, there also exists a reflection ordering $\varphi$ 
such that $\varphi ^{-1}(N)=s$ (cf.~(\ref{eq:ref-ord})).
Hence (\ref{eq:Rs=R}) implies that 
\[
\T s 
= \left( \prod _{i=1}^{N-1}  
(1+Q_{\varphi ^{-1}(i)}) \right)  (1+Q_s)\,  s
=\T\ . 
\]
It follows that, more generally,
$\T w = \T$
for all $w\in W$.
Analogously one shows that 
$w \T = \T$ for all $w\in W$.
These equations 
can be interpreted
as saying that the matrix of~$\T$ in the basis~$W$
of $\k[W]$ is invariant under permutations of rows and columns.
Hence there exists a constant~$c$ such that,
for any $u\in W$,
\[ \T (u) = c \sum _{v \in W} v  . \]
On the other hand, it is clear  from (\ref{eq:Qq}) 
that the coefficient of $\wnot$ in $\T(e)$ is~$\leq 1$,
where $e\in W$ is the identity element.
Since $c$ is obviously a positive integer, we conclude that $c=1$, 
and part~1 is proved. 
\medskip

To prove the rest, we will need the following lemma
that generalizes the corresponding result for the ordinary Bruhat
order (see, e.g.,~\cite[Lemma~4.1]{Dy}). 

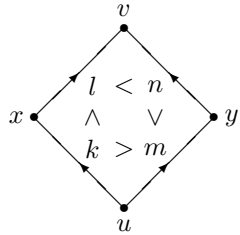
\begin{figure}[ht]
\setlength{\unitlength}{1.7pt} 
\begin{center}
\begin{picture}(40,45)(0,0)
\put(20, 0){\vector(1,1){10}}
\put(30, 10){\line(1,1){10}}
\put(20, 0){\vector(-1,1){10}}
\put(10, 10){\line(-1,1){10}}
\put(0, 20){\vector(1,1){10}}
\put(10, 30){\line(1,1){10}}
\put(40,20){\vector(-1,1){10}}
\put(30, 30){\line(-1,1){10}}

\put(0,20){\circle*{2}}
\put(20,0){\circle*{2}}
\put(20,40){\circle*{2}}
\put(40,20){\circle*{2}}

\put(20,-4){\makebox(0,0){$u$}}
\put(20,44){\makebox(0,0){$v$}}
\put(-4,20){\makebox(0,0){$x$}}
\put(44,20){\makebox(0,0){$y$}}

\put(13,13){\makebox(0,0){$k$}}
\put(13,27){\makebox(0,0){$l$}}
\put(27,13){\makebox(0,0){$m$}}
\put(27,27){\makebox(0,0){$n$}}

\put(13,20){\makebox(0,0){$\wedge$}}
\put(27,20){\makebox(0,0){$\vee$}}
\put(20,27){\makebox(0,0){$<$}}
\put(20,13){\makebox(0,0){$>$}}

\end{picture}
\end{center}
\caption{Tilted Bruhat interval of length 2}
\label{fig:diamond}
\end{figure}

\begin{lemma}
\label{lem:diamond} 
Assume that
\begin{equation}
\label{eq:a>b}
u,x,v \in W\ ,\quad u \stackrel{k}\to x \stackrel{l}\to v\ ,
\quad k>l\ .
\end{equation}
Then there exists $y \in W$ such that (cf.\ Figure~\ref{fig:diamond}) 
\begin{equation}
\label{eq:b<d>c<a}
u \stackrel{m}\to y \stackrel{n}\to v\ ,
\quad l<n>m<k\ .
\end{equation}
\end{lemma}

\proof 
Consider the dihedral group $W'$ generated by the reflections
$\varphi^{-1}(k)$ and $\varphi^{-1}(l)$. 
Define the operators $Q_\tau$ and $R_\tau$ 
by (\ref{eq:Qq})--(\ref{eq:Rq}), and write down the Yang-Baxter
equation~(\ref{eq:ybe}) for~$W'$, 
so that the order of the terms in the left-hand side was compatible
with the reflection order~$\varphi$. 
Thus the sequence of reflections appearing in the left-hand side 
is label-increasing,
while the one in the right-hand side is label-decreasing.
Apply the left-hand side to~$u$, and take the coefficient of~$\varepsilon^2 v$.
This will be the number of label-decreasing paths in $D(W)$ from $u$ to $v$ 
that have length~$2$ and stay within the coset~$W'u$.
We know one such path, namely 
$ u \stackrel{k}\to x \stackrel{l}\to v$.
By the Yang-Baxter equation, there should also be a label-increasing path
of length~$2$ from $u$ to $v$ that stays within~$W'u$;
let us denote it by $u \stackrel{m}\to y \stackrel{n}\to v$.
It remains to check that $k>m$ and $l<n$.
These two statements are completely analogous to each other,
so we will only show how to prove the first one.
Suppose that, on the contrary, $k<m$.
Then $l<k<m<n$, which in particular means that the four reflections
labelled by $l,k,m,n$ are all distinct.
If $W'$ is of type~$A_2\,$, this already brings the desired
contradiction,
since in that case there are only three reflections in~$W'$. 
If $W'$ is of type~$B_2\,$, with canonical generators $a$ and $b$
(say, $\varphi(a)<\varphi(b)$, then there are four reflections in
$W'$, and therefore $l,k,m,n$ correspond to $a,aba,bab,b$,
respectively. But this would imply that $v=a\cdot aba\cdot u = b\cdot
bab\cdot u$, a contradiction. 
The remaining case $W'=W=G_2$ is checked directly. 
\endproof

We can now complete the proof of Theorem~\ref{th:tilted} using an
argument borrowed from~\cite{Dy}.
Among all shortest paths in $D(W)$ from $u$ to $v$, let 
\begin{equation}
\label{eq:chain}
u=w_{1} \to w_{2} \to \cdots \to w_{d} =v
\end{equation} 
be the one whose label sequence is lexicographically minimal. 
To prove part~2 of the theorem, we need to show is that this path is
label-increasing.  
Suppose otherwise, i.e., for some $i \in \{ 2, \dots,d-1 \}$,
we have $w_{i-1}\stackrel{k} \to w_{i}\stackrel{l} \to w_{i+1}$
with $k>l$. 
(We cannot have $k=l$ since this would create a loop,
and the path would not be shortest.) 
Then, by Lemma~\ref{lem:diamond}, there exists $y\in W$ such that 
$w_{i-1}\stackrel{m} \to y \stackrel{n}\to w_{i+1}$ and 
$m <k$. 
Thus replacing $w_{i}$ by $y$ in~(\ref{eq:chain}) 
produces a chain with lexicographically smaller sequence 
of labels---a contradiction.  

Finally, part~3 follows from parts~1 and~2. 
\endproof

In the terminology of~\cite{BW2}, Theorem~\ref{th:tilted} asserts
that the tilted Bruhat order 
(hence any any tilted Bruhat interval $D(u,v)$) 
is \emph{EL-shellable} (hence \emph{CL-shellable}),
with the EL-shelling provided by any reflection ordering
(and therefore by its reversal as well). 

Theorem~\ref{th:tilted} 
implies the following generalization of Verma's theorem and its
refinements mentioned above.

\begin{corollary}
\label{cor:tilted}
Each tilted Bruhat order $D_u(W)$ of a Weyl group~$W$ is a
lexicographically shellable lower Eulerian poset. 

As a consequence, any tilted Bruhat interval is a face poset of a 
shellable regular CW sphere.
In particular, it is Eulerian and Cohen-Macaulay. 
\end{corollary}

\proof
By~\cite[Proposition~4.5]{bjorner-cw}, Theorem~\ref{th:tilted}
implies that $D(u,v)$ is a face poset of a regular CW sphere.
Such posets are well known to be both
Eulerian and Cohen-Macaulay; see, e.g.,
Stanley~\cite[Section~1]{StaNato}. 
\endproof

The Eulerian property can also be deduced directly from
Theorem~\ref{th:tilted} as follows. 
By a 
simple counting argument (cf.\ \cite[Corollary~2.3]{BGS}),
the values of the M\"obius function can be computed from an
EL-shelling by 
\[
\mu(u,v) = (-1)^{{\rm rank}(y)-{\rm rank}(x)}\,\cdot\,(\textrm{number
  of label-decreasing chains from $u$ to $v$}) \,.
\]
In our case, there is exactly one
such chain, and the Eulerian property follows.

\section*{Acknowledgments}

We thank Anders Bj\"orner, Alain Lascoux, and Richard Stanley 
for useful comments and conversations.
Part of our work was carried out while the first two authors were
participating in the ``Combinatorics'' program at MSRI.

\end{document}